\documentclass{article}
\usepackage[utf8]{inputenc}
\usepackage{epsf}
\usepackage{epsfig}
\usepackage{color}
\usepackage{amsmath}
\usepackage{amsthm}
\usepackage{authblk}
\usepackage{rotating}
\usepackage{lscape}
\usepackage[font=small,format=plain,labelfont=bf,up,textfont= up]{caption}
\usepackage{multirow}
\usepackage{graphicx}
\usepackage{subfigure}
\newtheorem{rmk}{Remark}[section]

\title{Membranes in Optic Nerve Models}
\author[1]{Yi Zhu}
\author[2,*]{Shixin Xu}
\author[3,4]{Robert S. Eisenberg}
\author[5,6,7,1]{Huaxiong Huang}

\affil[1]{Department of Mathematics and Statistics, York University, Toronto, ON, M3J 1P3, Canada}
\affil[2]{Duke Kunshan University, Jiangsu, Kunshan, 215316, China. }
\affil[3]{Department of Applied Mathematics, Illinois Institute of Technology, Chicago IL 60616 USA}
\affil[4]{Department of Physiology $\&$ Biophysics, Rush University Chicago IL 60612, USA}
 \affil[5]{Research Centre for Mathematics, Advanced Institute of Natural Sciences, Beijing Normal University (Zhuhai), 519087,  China}
\affil[6]{Division of Science and Technology, BNU- HKBU United International College, Zhuhai, 519087, China}%
\affil[7]{Department of Mathematics and Statistics, York University, Toronto, ON, M3J 1P3, Canada}
\affil[*]{Corresponding author: shixin.xu@dukekunshan.edu.cn}

\date{}

\begin{document}

\maketitle
\begin{abstract}
    Membranes are present in all cells and tissues. Mathematical models of cells and tissues need a compact mathematical description of membranes with a resolution of about 1 nm. Membranes isolate cells because ions have difficulty penetrating the dielectric barrier they create. Here we introduce a dielectric mathematical membrane condition to replace a condition that did not include dielectric properties. Our mathematical membrane condition includes a dielectric lipid bilayer punctured by channels that conduct ions selectively.
\end{abstract}
\section{Introduction}
Membranes isolate cells and organelles from the solutions around them. Membranes are lipid bilayers that isolate regions because of their polarization. The energy to move an ion from salt solution into the lipid of the bilayer is very large. As it moves, an ion must be separated from its ionic atmosphere (measured by the Debye length) and dielectric surround (measured by the Bjerrum length) that balance its permanent charge. The dielectric lipid has zero ionic atmosphere and little dielectric atmosphere to balance its charge. The movement of the ion from solution into the lipid of the membrane involves a huge energy barrier. The energy barrier isolates the interior of cells and organelles from their surround. 
Membranes are punctured by channels that provide pathways for ion flow with relatively small energy barriers (when the channel is open). In recent papers on the optic nerve \cite{zhu2021optic}, we followed Mori \cite{mori2011model,mori2015multidomain} and described the membranes by a  condition involving just the channels and their conductance

\begin{equation}
\boldsymbol{J_{k}^i}\cdot \boldsymbol{n} = \frac{g_m^i}{z^ie}(\phi_{in}-\phi_{ex}-E^i), 
\end{equation}
where $ i=\mathrm{Na}^{+},\mathrm{K}^{+},\mathrm{Cl}^{-} $, $E^i =\frac{k_BT}{z^ie} log\left(\frac{c_{ex}^i}{c_{in}^i}\right)$ is the Nernst potential. 

Here we introduce a mathematical membrane condition involving membrane dielectric properties (represented realistically as a capacitance) as well as membrane conductance. We show that our previous results are not significantly changed by introducing  the new mathematical membrane condition. We choose the resolution of the membrane condition so it describes membrane properties important in Hodgkin Huxley computations of the action potential and calculations of the associated ionic and water fluxes. Other effects, such as the ionic atmosphere surrounding the membrane, or the permanent charge inside the ion channel require higher resolution models.
We will use the realistic membrane condition (that includes capacitance) in future work lest readers and students misunderstand the physical basis of the main function of membranes, namely to provide a large dielectric energy barrier that isolates cells and organelles from their surrounds.

\section{Mathematical membrane condition}
In this section, we discuss the effect of the new mathematical membrane condition for both axon cell membrane and glial cell membrane.

The cell membrane can be viewed as a capacitor representing the lipid dielectric and a set of conductances representing ion channels following ample precedent.
In the model, we introduce the $J^{i}_{m,k}, (k=gl,ax, i=\mathrm{Na}^{+},\mathrm{K}^{+},\mathrm{Cl}^{-})$ for this ion flux and the transmembrane source term in the \cite{zhu2021optic} can be written as 
\begin{equation*}
J_{k}^{m, i}=J_{p,k}^{i}+J_{c,k}^{i}+J^{i}_{m,k}, \quad k=gl,ax, \quad i=\mathrm{Na}^{+}, \mathrm{K}^{+}, \mathrm{Cl}^{-}.
\end{equation*}
where 
\begin{equation*}
J^{i}_{m,k}= \lambda^{i}\frac{C_{m}}{z^{i}e} \frac{d (\phi_{k}-\phi_{ex})}{dt},
\end{equation*}
and active ion pump source  $J_{p,k}^{i}$  and  passive ion channel source $J_{c,k}^{i}$, on the $k$th membrane. The conservation of ion concentration implies the following system of partial differential equations to describe the dynamics of ions in each region, for $i=\mathrm{Na^+,K^+,Cl^-}$  
\begin{equation*}
\label{ion_governing}
\begin{aligned}
&\frac{\partial\left(\eta_{g l} c_{gl}^{i}\right)}{\partial t}+\mathcal{M}_{gl}J^{m,i}_{gl}+\nabla \cdot\left(\eta_{g l} \boldsymbol{j}_{gl}^{i}\right)=0, ~~~&&\text { in } \Omega_{OP}, \\
&\frac{\partial\left(\eta_{ax} c_{ax}^{i}\right)}{\partial t}+\mathcal{M}_{ax}J^{m,i}_{ax}+\frac{\partial}{\partial z}\left(\eta_{ax} j_{ax,z}^{i}\right)=0,  ~~~~~&&\text { in } \Omega_{OP}, \\
&\frac{\partial\left(\eta_{ex} c_{ex}^{i}\right)}{\partial t}-\mathcal{M}_{ax}J^{m,i}_{ax}-\mathcal{M}_{gl}J^{m,i}_{gl}+\nabla \cdot\left(\eta_{ex} \boldsymbol{j}_{e x}^{i}\right)=0, ~~~~ &&\text { in } \ \Omega_{OP},
\end{aligned}
\end{equation*}
The electric potentials in the axon compartment, glial compartment and extracellular space are 
\begin{equation*}
\label{phi_governing}
\begin{aligned}
\sum_{i} z^{i} e \mathcal{M}_{gl}J^{m,i}_{gl}+\sum_{i} z^{i} e \nabla \cdot\left(\eta_{g} \boldsymbol{j}_{gl}^{i}\right) &=0, \\
\sum_{i} z^{i} e \mathcal{M}_{a x}J^{m,i}_{ax}+\sum_{i} z^{i} e \frac{\partial}{\partial z}\left(\eta_{ax} j_{ax, z}^{i}\right) &=0, \\
\sum_{i} z^{i} e \nabla \cdot\left(\eta_{gl} \boldsymbol{j}_{gl}^{i}\right)+\sum_{i} z^{i} e \frac{\partial}{\partial z}\left(\eta_{a x} j_{ax,z}^{i}\right)+\sum_{i} z^{i} e \nabla \cdot\left(\eta_{ex} \boldsymbol{j}_{ex}^{i}\right) &=0,
\end{aligned}
\end{equation*}

\section{Simulation}
Here we repeat the simulations of the optic nerve published previously but with modified parameters.
We choose the capacitance on both axon and glial membrane as $C_{m}=7.5 \times 10^{-3} \ \mathrm{F\cdot m^{-2}}$ \cite{major1994detailed} and $\lambda^{i}=\frac{1}{3}$ for each ion specie. Furthermore, action potentials are generated by the axon with the parameters in the following Table  \ref{tab:my_label}.

\begin{table}[hpt]
    \centering
    \begin{tabular}{|c|c|c|}
    \hline
       Parameters   & Previous Value  &  New Value   \\
     \hline
        $\mathcal{M}_{ax}$ &  $5.98 \times 10^{6} \mathrm{m}^{-1}$  & $2.392\times 10^{5} \mathrm{m}^{-1}$  \\
      \hline
       $I_{ax,1}$ & $9.56 \times 10^{-4} \mathrm{~A} / \mathrm{m}^{2}$  &  $ 2.39\times 10^{-2} \mathrm{~A} / \mathrm{m}^{2}$ \\
       \hline 
       $I_{ax,2}$ & $1.3 \times 10^{-4} \mathrm{~A} / \mathrm{m}^{2}$  & $3.25 \times 10^{-3} \mathrm{~A} / \mathrm{m}^{2}$ \\ 
    \hline
    $g_{leak}^{Na}$ & $4.8 \times 10^{-3} \mathrm{~S} /  \mathrm{m}^{2}$ & $1.2\times 10^{-1} \mathrm{~S} /  \mathrm{m}^{2}$  \\
     \hline
    $g_{leak}^{K}$ & $2.2 \times 10^{-2} \mathrm{~S} / \mathrm{m}^{2}$   &  $5.5\times 10^{-1} \mathrm{~S} /  \mathrm{m}^{2}$ \\ 
    \hline
    $\bar{g}^{Na}$ & $1.357 \times 10^{1} \mathrm{~S} / \mathrm{m}^{2}$ & $3.393 \times 10^{2} \mathrm{~S} / \mathrm{m}^{2}$ \\
    \hline 
    $\bar{g}^{K}$ & $2.945 \mathrm{~S} / \mathrm{m}^{2}$ & $7.364 \times 10^{1} \mathrm{~S} / \mathrm{m}^{2}$ \\
    \hline
    $g_{ax}^{Cl}$ & $1.5 \times 10^{-1} \mathrm{~S} / \mathrm{m}^{2}$ & $3.75 \mathrm{~S} / \mathrm{m}^{2}$ \\
    \hline
    $I_{shock}$ & $3\times 10^{-3} \mathrm{~A\cdot m^{-2}}$  &  $7.5 \times 10^{-2} \mathrm{~A\cdot m^{-2}}$\\
    \hline
    \end{tabular}
    \caption{Parameters Table}
    \label{tab:my_label}
\end{table}
\begin{rmk}
 
The parameters are adjusted based on the following facts
\begin{enumerate}
    \item The size of necturus optic nerve is much smaller that squid;
    \item In \cite{zhu2021optic}, we estimated the membrane conductance based on the experimental results \cite{orkand1966effect} and  the fact that the density of ion channel in Necturus ($24/\mu m^2$) is lower that that of squid ($170-290/\mu m^2$). Without the capacitance effect, the conductance of axon may be underestimated.   
\end{enumerate}
At the same time, in \cite{orkand1966effect}, the  measured $\mathrm{K^+}$ variation is around $0.2mM$ for single action potential. In order to match this, we keep the value $\mathcal{M}_{ax}*\bar{g}_k$ to be same, while decrease the value of $\mathcal{M}_{ax}$ and increase the value of conductance. The other parameters are adjusted based the rest state ion concentrations.  
\end{rmk}
\begin{figure}[h]
	\centering
   \includegraphics[width=5in,height=7cm]{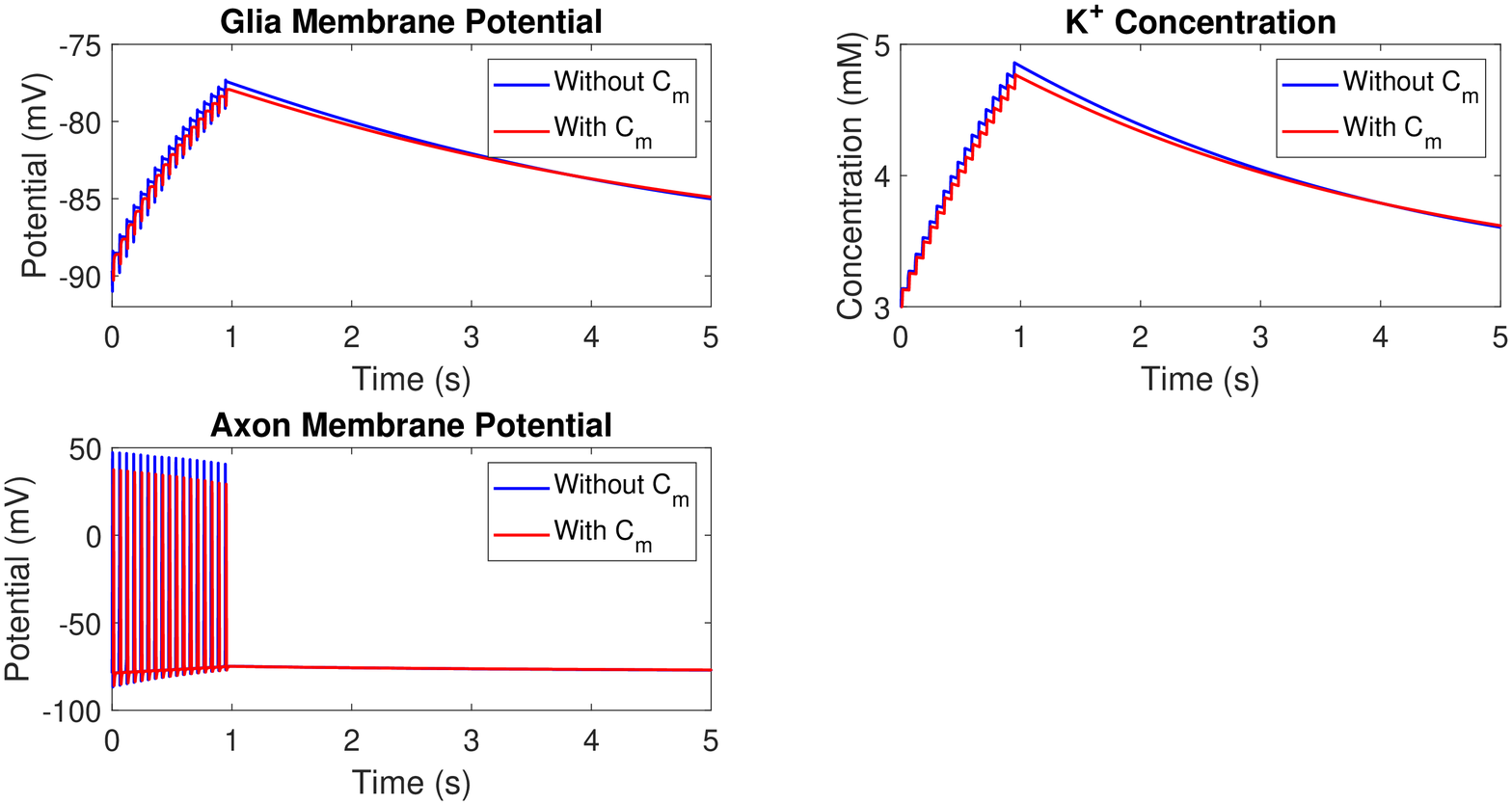}
	\caption{Comparison between the model with and without cell membrane capacitance effect.}
	\label{fig::comparision}
\end{figure}
Fig. \ref{fig::comparision} shows results from simulations using  the parameters in Table 1. We simulate  Orkand’s experimental setup with the bath solution containing  $\mathrm{3mM~ K}^+$. In Figure 1, we  compare the glial membrane potential, potassium concentration in the extracellular space and axon membrane potential.  The results do not differ significantly from those published using Mori’s purely conductive model of the membrane.


\end{document}